\documentclass[fleqn]{amsart}
\usepackage{amssymb, graphics, amsmath, epsfig, changebar, hyperref}
\usepackage{enumitem}
\usepackage{mathtools} 
\usepackage[percent]{overpic}
\usepackage[caption=false]{subfig} 

\theoremstyle{plain}
\newtheorem{theorem}{Theorem}
\newtheorem{corollary}{Corollary}

\theoremstyle{definition}

\theoremstyle{example}

\theoremstyle{remark}

\numberwithin{equation}{section}
\newcommand{\N}{\mathbb{N}}

\newcommand{\C}{\mathbb{C}}

\newcommand{\Om}{\Omega}

\graphicspath{{./images/}}

\begin{document}                   

\title{On two crossing numbers of algebraic knots under Hopf fibration}
\author{Maciej Mroczkowski}
\address{Institute of Mathematics\\
Faculty of Mathematics, Physics and Informatics\\
University of Gdansk, 80-308 Gdansk, Poland\\
e-mail: mmroczko@mat.ug.edu.pl}

\begin{abstract}
We answer a question posed by Fielder in \cite{F} concerning two notions of crossing number for algebraic knots $K$ under Hopf fibration, one topological, denoted $h(K)$, the other coming from the realization of such knots around complex singularities, denoted $C_{alg}(K)$. We show that $C_{alg}(K)-h(K)$ can be arbitrarily large.
We also give an upper bound for $h$ of some families of knots such as torus knots $T(2,n)$, twist knots and their mirror images.
\end{abstract}
\maketitle

\let\thefootnote\relax 

\section{Introduction}

Let $p:S^3\to S^2$ be the Hopf fibration.
Given a link $L$ in $S^3$, let $h(L)$ be the minimal number of crossings of
$p(L)$ among all generic projections of $L$ with $p$. We call $h(L)$ the
{\it Hopf crossing number} of $L$. It was first considered by Fiedler in
\cite{F} and studied by the author in \cite{M1,M2}.

For algebraic links $L$ there is a number defined analogously to $h$, coming
from the {\it canonical Hopf fibration}, denoted $C_{alg}(L)$ and studied by
Fiedler in \cite{F}. Let the algebraic link $L$ be realized as the 
intersection of a complex plane algebraic curve $X$ and a small sphere $S^3$
with a singularity of $X$ at its center. The canonical Hopf fibration
from $S^3$ to $\C P^1$ is the fibration induced by the intersection of $S^3$
and the complex lines passing through the center of $S^3$. Fiedler defined
$C_{alg}(L)$ as the minimal number of crossings under such a canonical Hopf
fibration among all generic algebraic realizations of $L$.

Recall, that an algebraic knot of type $\{(p_1,q_1),(p_2,q_2)\}$ is a $(p_2,q_2)$ cable of the torus knot of type $(p_1,q_1)$. Also, for $i\in\{1,2\}$, $1<p_i<q_i$, $gcd(p_i,q_i)=1$ and $q_2>p_1p_2q_1$ (see\cite{Tr}).
Fiedler has shown that the following inequality holds for such algebraic knots:

\begin{theorem}[Theorem 2 in \cite{F}]\label{thm_fiedler}
Let $K$ be an algebraic knot of type $\{(p_1,q_1),(p_2,q_2)\}$. If $q_1<2p_1$ and both $p_1$ and $q_1$ are odd, then:
\[C_{alg}(K)\ge q_2-q_1 p_2 p_1 \]
\end{theorem}

It is clear that for any algebraic link $h(L)\le C_{alg}(L)$. In the paragraph following this theorem, Fiedler posed the following problem, which we rephrase (Fiedler did not use the notation $h$): For $K$ as in the preceding  theorem, is it true that $h(K)<
q_2-q_1 p_2 p_1$? More generally, is always $h(K)=C_{alg}(K)$?

We answer the first question in Theorem~\ref{thm:main} for $p_2=2$ and $q_2$ sufficiently large (larger than a constant depending on $p_1$ and $q_1$). Under these assumptions $h(K)<q_2-2 p_1 q_1$. This gives obviously a negative answer to the second question: for infinite families of algebraic knots we have $h(K)<C_{alg}(K)$. In fact, we show that the gap between $h(K)$ and $C_{alg}(K)$ can be arbitrarily large. 

In section~\ref{sec:arr} we briefly recall the definition of arrow diagrams and their relationship to the Hopf crossing number.
In section~\ref{sec:upper} we give upper bounds for the Hopf crossing number of knots (or links) obtained from a $4$-tangle by closing it with some half-twists (Theorem~\ref{thm:h_twist}). In particular, we get upper bounds for $h$ of the torus knots $T(2,n)$, the twist knots and their mirror images. Finally, we answer Fiedler's question (Theorem~\ref{thm:main}).

\section{Arrow diagrams of links in $S^3$}\label{sec:arr}
\begin{figure}[h]
\centering
\includegraphics{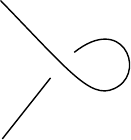} \raisebox{2 em}{$\; \underset{\longleftrightarrow}{\Omega_1} \;$} \includegraphics{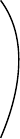}
\hspace{2em}
\includegraphics{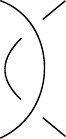} \raisebox{2 em}{$\; \underset{\longleftrightarrow}{\Omega_2} \;$} \includegraphics{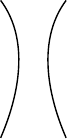}
\hspace{2em}
\includegraphics{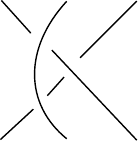} \raisebox{2 em}{$\; \underset{\longleftrightarrow}{\Omega_3} \;$} \includegraphics{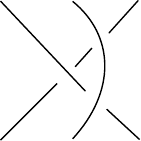}
\hspace{2em} \\[2em]
\includegraphics{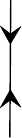} \raisebox{2 em}{$\; \underset{\longleftrightarrow}{\Omega_4} \;$} \includegraphics{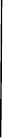}
\raisebox{2 em}{$\; \underset{\longleftrightarrow}{\Omega_4} \;$} \includegraphics{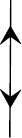} \hspace{5em}
\includegraphics{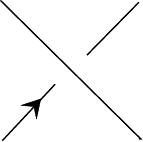} \raisebox{2 em}{$\; \underset{\longleftrightarrow}{\Omega_5} \;$} \includegraphics{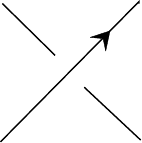}
\\[2em]
\includegraphics{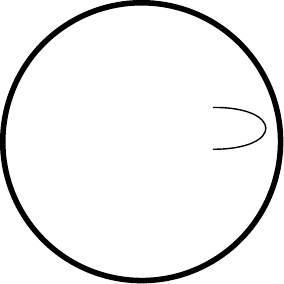} \raisebox{4.2 em}{$ \;
 \underset{\longleftrightarrow}{\Omega_{\infty}} \; $ } \includegraphics{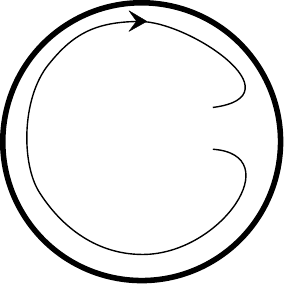}
\caption{Reidemeister moves}
\label{reid_moves}
\end{figure}

It was shown in \cite{M1} that projections of knots (or links) under the Hopf fibration $p:S^3\to S^2$ can be
encoded as arrow diagrams in a disk, two such diagrams representing the same knot if and only if one can be transformed into the other with a series of six Reidemeister moves shown in Figure~\ref{reid_moves}. An arrow diagram looks like a
classical diagram except that there may be some arrows on its strands outside the crossings. These arrows are involved
in Reidemeister moves $\Om_4$, $\Om_5$ and $\Om_\infty$ (for this last move the boundary of the disk is drawn in thick).
See \cite{MD} for the interpretation of the arrow diagrams and Reidemeister moves in a more general context (links in a product
of a surface times $S^1$).

For a link $L$ in $S^3$, the Hopf crossing number of $L$, or $h(L)$, is the minimum number of crossings among all arrow
diagrams of $L$. Since classical diagrams are also arrow diagrams with no arrows, it is obvious that $h(L)\le c(L)$, where
$c(L)$ is the classical crossing number of $L$.

\section{Upper bounds for the Hopf crossing number}\label{sec:upper}
\begin{figure}[h]
\centering
\includegraphics{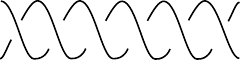}
\caption{$6$ positive half-twists}
\label{twist6}
\end{figure}

By $n$ positive half-twists, we mean a four tangle, such as the one shown in Figure~\ref{twist6}, for $n=6$. 
By switching all crossings, we get $n$ negative half-twists.

\begin{figure}[h]
\centering
\includegraphics{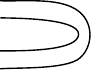}\hspace{3em}\includegraphics{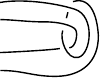}\hspace{3em}\includegraphics{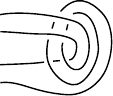}
\caption{Tangles $G_0$,  $G_1$ and $G_2$}
\label{t1_2}
\end{figure}

\begin{figure}[h]
\centering
\includegraphics{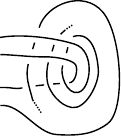} \raisebox{1.9 em}{$\; \longrightarrow \;$}
\includegraphics{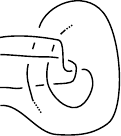}
\raisebox{1.9 em}{$\; \longrightarrow \;$}
\includegraphics{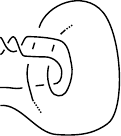}
\caption{$G_n$ to $G_{n-1}$}
\label{tg_n}
\end{figure}

Let us define a family of $4$-tangles, $G_n$, $n\in\N$, as in Figure~\ref{t1_2}, for $n=0,1,2$ and
on the left of Figure~\ref{tg_n} in general. In this last figure, it is shown that, for $n>0$,
$G_n$ can be transformed into $G_{n-1}$ with an additional positive twist (two positive half-twists).
Thus, $G_n$ can be transformed into $G_0$ with $n$ postive twists ($2n$ positive half-twists).

\begin{theorem}\label{thm:h_twist}
Let $K$ be a knot (or link) obtained from a tangle $T$ with $c$ crossings, by adding $2n+1$ positive half-twists, $n\ge 2$.
Then $h(K)\le c+n-1$.

If the $2n+1$ half-twists are negative, then $h(K)\le c+n$.

\begin{proof}
Assume, for example, that $n=3$. Such $K$ is shown in Figure~\ref{knotK}.

\begin{figure}[h]
\centering
\includegraphics{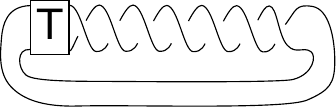}
\caption{knot $K$ obtained from $T$ and $7$ positive half-twists}
\label{knotK}
\end{figure}

Insert the tangle $T$ into a diagram with two arrows as on the left of Figure~\ref{t27}. In this figure, the diagram is transformed
in four steps:

\begin{enumerate}[label=(\roman*)]
\item The exterior arrow is eliminated: create a kink with $\Om_1$, push the arrow inside the kink with $\Om_5$, use
$\Om_\infty$ and $\Om_4$ to eliminate the arrow.
\item The second arrow is pushed outside and eliminated: use three $\Om_2$ moves and three $\Om_5$ moves to push
the arrow outside. Then eliminate the arrow as in the previous step. Now we have an extra arc going around the whole diagram.
Since there are no more arrows it can be moved above the rest of the diagram with classical Reidemeister moves, then removed
with $\Om_1$.
\item The tangle $T$ glides in a counterclockwise direction on two strands, rotating at the same time with a full twist.
\item An arc is moved above some part of the diagram.
\end{enumerate}

\begin{figure}[h]
\centering
\includegraphics{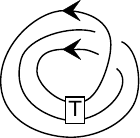} \raisebox{1.9 em}{$\; \overset{i}{\rightarrow} \;$} \includegraphics{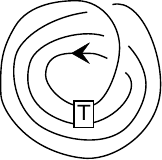}
\raisebox{1.9 em}{$\; \overset{ii}{\rightarrow} \;$} \includegraphics{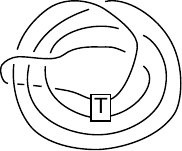}
\raisebox{1.9 em}{$\;\overset{iii}{\rightarrow} \;$} \includegraphics{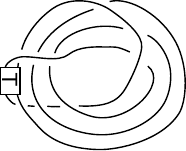}
\raisebox{1.9 em}{$\;\overset{iv}{\rightarrow} \;$} \includegraphics{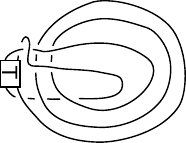}
\caption{Eliminating arrows, moving tangle}
\label{t27}
\end{figure}

At the end we get an arrowless diagram with $T$, $3$ positive half-twists and the tangle $G_2$. As was noticed before this theorem, $G_2$ can be transformed into $4$ positive half-twists. Thus, we get a diagram of $T$ with $7$ positive half-twists. The initial arrow diagram has $c+2$ crossings. Thus, $h(K)\le c+2$. 

The same argument works for any $n\ge 2$: start with an arrow diagram such as the one on the left of Figure~\ref{t27}, with the curve going around the center $n$ times (in the figure $n=3$), the two arrows located on the outermost and innermost loops and $T$ on the innermost loop and the next innermost one.  Such an arrow diagram has $c+n-1$ crossings. As before, we transform it into $T$, $3$ positive half-twists and the tangle $G_{n-1}$. This, in turn can be transformed into $T$ with $2n+1$ positive half-twists.. Hence, $h(K)\le c+n-1$.

If the $2n+1$ half-twists in $K$ are negative the situation is similar. The difference is that we start from a diagram obtained from the second diagram in Figure~\ref{t27} by switching all crossings except those in $T$ and changing the direction of the interior arrow.
Such a diagram has $c+n$ crossings.
Removing the arrow is similar to step $ii$: one gets the third diagram in the figure with all crossings (outside $T$) switched. The next two steps
are mirroring steps $iii$ and $iv$. We get an arrowless diagram with $T$, $3$ negative half-twists and the mirror image of $G_n$, denoted $\overline{G_n}$.
By mirroring Figure~\ref{tg_n}, $\overline{G_n}$ can be transformed into $G_0$ with $2n$ negative half-twists. 
Thus we get $T$ with $2n+1$ negative half-twists, from an arrow diagram with $c+n$ crossings. Hence, $h(K)\le c+n$.
\end{proof}
\end{theorem}

\begin{corollary}\label{cor:torustwist}
Let $T(2,n)$, $n\ge 5$ be the right-handed $(2,n)$ torus knot and $\overline{T(2,n)}$ its mirror image.
Let $K_n$, $n\ge 7$, be the twist knot and $\overline{K_n}$ its mirror image.
Then we have:
\[h(T(2,n))\le \frac{n-3}{2},\quad h(\overline{T(2,n)})\le \frac{n-1}{2},\]
\[h(K_n)\le \left\lfloor\frac{n+1}{2}\right\rfloor,\quad h(\overline{K_n})\le \left\lfloor\frac{n+1}{2}\right\rfloor-1\]
\begin{proof}
For $T(2,n)$ (resp. $\overline{T(2,n)}$), we use Theorem~\ref{thm:h_twist} with $T$ a tangle with $0$ crossings and $n$ positive (resp. negative) half-twists. As $n$ is odd, $n=2\left(\frac{n-1}{2}\right)+1$, and we get $h(T(2,n))\le \frac{n-1}{2}-1=\frac{n-3}{2}$ and $h(\overline{T(2,n)})\le \frac{n-1}{2}$.

The knot $K_n$ (resp. $\overline{K_n}$) is obtained from a tangle with $2$ crossings adding $n-2$ negative (resp. positive) half-twists. Now, $n-2=2\left\lfloor \frac{n-3}{2} \right\rfloor+i$ with $i=1$ if $n$ is odd and $i=2$ if $n$ is even. From Theorem~\ref{thm:h_twist}, it follows that $h(K_n)\le 2+\left\lfloor \frac{n-3}{2} \right\rfloor=\left\lfloor \frac{n+1}{2} \right\rfloor$ and, similarly, $h(\overline{K_n})\le \left\lfloor\frac{n+1}{2}\right\rfloor-1$.
\end{proof}
\end{corollary}

In fact, using results from \cite{M2}, it can be shown that $h(T(2,5))=1$, $h(\overline{T(2,5)})=2$,  $h(T(2,7))=2$, $h(\overline{T(2,7)})=3$, $h(K_7)=3$ and $h(\overline{K_7})=2$.
From this small set of values of $h$ it looks like the inequalities in the last corollary may be best possible for the torus knots, but not for the twist knots. 

It was noticed in Section~\ref{sec:arr}, that $h(K)\le c(K)$ for any knot (or link) $K$, where $c(K)$ is the classical crossing number. It is interesting to compare these two invariants, by studying $c(K)-h(K)$ and $\frac{h(K)}{c(K)}$. For example, for $(n,n+1)$ torus knots,  $h(T(n,n+1))=0$ (see~\cite{M1}), so, obviously, $\frac{h(T(n,n+1))}{c(T(n,n+1))}=0$.
Since $c(T(2,n))=c(K_n)=n$ (and the same for their mirror images), we get an immediate corollary:

\begin{corollary}
Let $L_n$ be any of the four families of knots: $T(2,n)$, $\overline{T(2,n)}$, $K_n$ or $\overline{K_n}$. Then:
\[\lim_{n\to\infty} c(L_n)-h(L_n)=\infty, \quad
\limsup_{n\to\infty}\frac{h(L_n)}{c(L_n)}\le\frac{1}{2} \]
\end{corollary} 

Notice, in particular, that these are all alternating knots. Other such examples of alternating knots could be easily constructed by taking other tangles $T$ in Theorem~\ref{thm:h_twist}, so that adding half-twists, we obtain alternating knots. For all such families of alternating knots, there is no bound on $c(K)-h(K)$. 

We ask the question: Is there a family of alternating knots $L_n$ and some constant $k<\frac{1}{2}$ satisfying $\displaystyle\lim_{n\to\infty} c(L_n)=\infty$ and $\displaystyle\limsup_{n\to\infty}\frac{h(L_n)}{c(L_n)}\le k$ ? What is the lowest such $k$ possible?

Finally, we state our main theorem, answering Fiedler's question.

\begin{theorem}\label{thm:main}
Let $K$ be an algebraic knot of type $\{(p_1,q_1),(2,q_2)\}$ with $p_1$ and $q_1$ odd (and $q_2$ odd since $gcd(2,q_2)=1$).
Let $Q_2=10p_1q_1-6q_1-1$. Then for any $q_2\ge Q_2$ (in particular $q_2>2p_1q_1$ so $K$ is algebraic), one has: 
\[h(K)<q_2-2p_1q_1\le C_{alg}(K)\]
In fact, the gap between $C_{alg}(K)$ and $h(K)$ can be arbtrarily large since, for $q_2\ge Q_2$, $C_{alg}(K)-h(K)\ge \frac{q_2-Q_2}{2}+1$.

In particular $Q_2=119$ for $(p_1,q_1)=(3,5)$, so $h(K)<C_{alg}(K)$ for the algebraic knot of type $\{(3,5),(2,119)\}$. Also, for $n\ge 1$, $C_{alg}(K)-h(K)\ge n$ for the algebraic knot of type $\{(3,5),(2,117+2n)\}$.
\begin{proof}
Taking $p_2=2$ in Theorem~\ref{thm_fiedler}, we get $C_{alg}(K)\ge q_2-2p_1q_1$.
The torus knot of type $(p_1,q_1)$ has a diagram with $(p_1-1)q_1$ crossings. For $q_2>2p_1q_1$, a $(2,q_2)$ cabling of such a knot has a diagram with $4(p_1-1)q_1+q_2-2(p_1-1)q_1$
crossings. The term with minus comes from the fact that, to get the $(2,q_2)$ cable, one starts with a diagram of the torus knot which has $(p_1-1)q_1$ positive crossings, so one adds as many negative
kinks to get a diagram with $0$ framing. Each such negative kink gives rise to two negative half-twists in the cable which cancel with some of the $q_2$ positive half-twists.

Thus $K$ consists of a tangle with $4(p_1-1)q_1$ crossings and $q_2-2(p_1-1)q_1$ positive half-twists.
From Theorem~\ref{thm:h_twist}, $h(K)\le 4(p_1-1)q_1+\frac{q_2}{2}-(p_1-1)q_1-\frac{3}{2}$.
Thus, $h(K)\le \frac{q_2}{2}+3(p_1-1)q_1-\frac{3}{2}$.

It follows that, as functions of $q_2$, $C_{alg}(K)$ is bounded from below by $f_1(q_2)=q_2-C_1$ and $h(K)$ is bounded from above by $f_2(q_2)=\frac{q_2}{2}+C_2$ with positive constants $C_1=2p_1q_1$ and $C_2=3(p_1-1)q_1-\frac{3}{2}$. The linear functions $f_1$ and $f_2$ are equal in $q_2=2(C_1+C_2)$.
Let $Q_2=2(C_1+C_2)+2=10p_1q_1-6q_1-1$. Then $q_2\ge Q_2$ implies $C_{alg}(K)\ge q_2-2p_1 q_1>h(K)$.
Furthermore $f_1(q_2)-f_2(q_2)\ge \frac{q_2-Q_2}{2}+1$, so $C_{alg}(K)-h(K)\ge \frac{q_2-Q_2}{2}+1$.
\end{proof}
\end{theorem}

We remark, that for $K$ as in the previous theorem,  neither $C_{alg}(K)$ nor $h(K)$ are known exactly. It may be true that $C_{alg}(K)>h(K)$ for some $q_2$ that are lower than $Q_2$. In fact, even for torus knots, it is only known that $C_{alg}(K)=h(K)=k-1$, when $K=T(n,n+k)$ for $k\le 3$, with the exception $C_{alg}(T(2,5))=h(T(2,5))=1$ (see \cite{F},\cite{M1}).


\begin{thebibliography}{99}
\bibitem{F} T. Fiedler:
{\it Algebraic links and the Hopf Fibration},
Topology {\bf 30} (1991), nr. 2, 259-265.


\bibitem{MD} M. Mroczkowski and M. Dabkowski, {\it KBSM of the product 
of a disk with two holes and $S^1$}, Topology and its Applications {\bf 
156} (2009), 1831--1849.

\bibitem{M1} M. Mroczkowski, {\it Knots with Hopf crossing number at most one},
Osaka J. Math. {\bf 57} (2020), no. 2, 279--304,
 \href{https://arxiv.org/abs/1907.11614}{arXiv:1907.11614}.

\bibitem{M2} M. Mroczkowski, {\it On some moves on links and the Hopf crossing number}, 
\href{https://arxiv.org/abs/1908.00342}{arXiv:1908.00342}.
accepted in Mediterranean Journal of Mathematics.


\bibitem{Tr} Le Dung Trang, {\it Sur les noeuds algebriques}, Compositio Mathematica {\bf 25} (1972): 281--321.



\end{thebibliography}
\end{document}